\documentclass[a4paper,twoside,11pt]{article}
\usepackage{bbm}
\usepackage{amsfonts}
\usepackage{amssymb}\usepackage{amsmath}
 \numberwithin{equation}{section}
\begin{document}

\title{Blow-up phenomena for the rotation-two-component Camassa-Holm system}

\author{Jingjing Liu\footnote{e-mail:
jingjing830306@163.com }\\
Department of Mathematics and Information Science, \\
Zhengzhou University of Light Industry,\\
 450002 Zhengzhou, China}
\date{}
\maketitle

\begin{abstract}
In this paper, we study blow-up phenomena for the rotation-two-component
Camassa-Holm (R2CH) system, a model of equatorial water waves that includes the effect of the
Coriolis force. We first derive the precise blow-up scenario for R2CH with $\sigma>0.$  Then, we present two new
blow-up results for strong solutions to the system and then the precise blow-up rate
for strong solutions to the system is given.\\

\end{abstract}

\section{Introduction}
\ \ \ \ \ \ In this paper, we consider the Cauchy problem of the following
rotation-two-component Camassa-Holm (R2CH) system, which was recently derived by Fan, Gao and Liu\cite{f-g-l}:
\begin{equation}\label{eq:original}
\left\{\begin{array}{ll}u_{t}-u_{txx}-Au_{x}+3uu_{x}=\sigma(2u_{x}u_{xx}+uu_{xxx})-\mu u_{xxx}\\
\ \ \ \ \ \ \ \ \ \ \ \ \ \ \ \ \ \ \ \ \ \ \ \ \ \ \ \ \ \ \ \ \ \ \ \ \ \ \ \ \ \ \ \ \ \ \ -(1-2\Omega A)\rho\rho_{x}+2\Omega\rho(\rho u)_{x},&t > 0,\,x\in \mathbb{R},\\
 \rho_{t}+(u\rho )_{x}=0, &t > 0,\,x\in \mathbb{R},\\
u(0,x) =u_{0}(x),&x\in \mathbb{R}, \\
\rho(0,x) = \rho_{0}(x),&x\in \mathbb{R},\end{array}\right.
\end{equation}
where the function $u(t,x)$ is the fluid velocity in the
$x$-direction, $\rho(t,x)$ is related to the free surface elevation from equilibrium, the parameter $A$ characterizes a linear underlying shear flow, the real dimensionless constant $\sigma$ is a parameter which provides the competition, or balance, in fluid convection between nonlinear steepening and amplification due to stretching, $\mu$ is a nondimensional parameter and $\Omega$ characterizes the constant rotational speed of the Earth. The boundary assumptions associated with (1.1) are $u\rightarrow0,$ $\rho\rightarrow1$ as $|x|\rightarrow\infty.$ The system (1.1) models the equatorial water waves with the effect of the Coriolis force. It is known that the gravity and the Coriolis force induced by the Earth¡¯s rotation are primary influence forces for the geophysical water waves. However, the effect of the Coriolis force is small since the smallness of the variation in latitude of the Equatorial Undercurrent in the equatorial ocean regions. So, the authors in \cite{f-g-l} use the approximation by the $f$-plane governing equations.

In the absence of the Earth's rotation, i.e., $\Omega=0,$ the system (1.1) becomes the generalized Dullin-Gottwald-Holm system \cite{c-l, h-g-g}
\begin{equation}\label{eq:original}
\left\{\begin{array}{ll}u_{t}-u_{txx}-Au_{x}+3uu_{x}=\sigma(2u_{x}u_{xx}+uu_{xxx})-\mu u_{xxx}-\rho\rho_{x},&t > 0,\,x\in \mathbb{R},\\
 \rho_{t}+(u\rho )_{x}=0, &t > 0,\,x\in \mathbb{R},\\
\end{array}\right.
\end{equation}
which was derived in the shallow-water regime following Ivanov's modeling approach \cite{i}. If we further take $\sigma=1$
in (1.2), it then becomes the standard two-component Dullin-Gottwald-Holm system
\begin{equation}\label{eq:original}
\left\{\begin{array}{ll}u_{t}-u_{txx}-Au_{x}+3uu_{x}=2u_{x}u_{xx}+uu_{xxx}-\mu u_{xxx}-\rho\rho_{x},&t > 0,\,x\in \mathbb{R},\\
 \rho_{t}+(u\rho )_{x}=0, &t > 0,\,x\in \mathbb{R}.\\
\end{array}\right.
\end{equation}
Recently, the generalized Dullin-Gottwald-Holm system (1.2) has been
studied in \cite{cy,c-l,h-g-g}. These works established local well-posedness on the line \cite{c-l,h-g-g}
and on the circle \cite{cy} respectively. In \cite{c-l}, the authors showed that (1.2) can still exhibit the wave-breaking
phenomenon, determined the exact blow-up rate of such solutions and established a sufficient condition for global solutions. In \cite{h-g-g},
the authors classified all traveling wave solution of this system, studied the blow-up mechanism and given two sufficient conditions which
can guarantee wave-breaking phenomena. For the periodic case \cite{cy}, the precise blow-up scenarios
of strong solutions and several results of blow-up solutions with certain initial profiles
are described in detail, the exact blow-up rates are also determined and a
sufficient condition for global solutions is established.

For $\mu=0$ in (1.3), it becomes the famous two-component Camassa-Holm
\begin{equation}\label{eq:original}
\left\{\begin{array}{ll}u_{t}-u_{txx}-Au_{x}+3uu_{x}=2u_{x}u_{xx}+uu_{xxx}-\rho\rho_{x},&t > 0,\,x\in \mathbb{R},\\
 \rho_{t}+(u\rho )_{x}=0, &t > 0,\,x\in \mathbb{R}£¬\\
\end{array}\right.
\end{equation}
witch was proposed first by Chen et al in \cite{clz} and Falqui in \cite{f}. Then, Constantin and Ivanov derived it in \cite{ci} in the
context of shallow water regime. Notice that the Camassa-Holm equation \cite{ch} can be obtained via the obvious reduction $\rho=0$ and $A=0.$
The Cauchy problem of (1.4) has been extensively studied \cite{c, ce1,ce2,ce3,ekl,ely,gy,gl1,gl2,gp,wwz,yl,yin}, here, we only list some of the results. It has been shown this system is locally
well-posed with initial data $(u_{0},\rho_{0})\in H^{s}\times H^{s-1},$ $s\geq2$  by Kato¡¯s semigroup theory \cite{ely}. Then the authors in \cite{gl1}
improved this result with initial data in the Besov spaces (specifically $(u_{0},\rho_{0})\in H^{s}\times H^{s-1},$ $s\geq\frac{3}{2}$ ). More interestingly,
it has global strong solutions \cite{gy, gl2} and also finite time blow-up solutions \cite{ely, gy, gl1}.

The system (1.1) has been studied on the line \cite{cfg,f-g-l,  zl}
and on the circle \cite{m, z} respectively. These works established the local well-posedness in
$H^{s}\times H^{s-1},$ $s> \frac{3}{2},$ derived wave-breaking scenario: the solution
blows up at time $T$ if and only if
$$\lim_{t\rightarrow T^{-}}\left\{\sup_{x\in \mathbb{R}}|u_{x}(t,x)|\right\}=+\infty.$$
and
\begin{equation}\label{eq:original}
\lim_{t\rightarrow T^{-}}\left\{\inf_{x\in \mathbb{R}}u_{x}(t,x)\right\}=-\infty, \ \ \ \ \text{for} \ \ \ \sigma=1, \mu=0
\end{equation}
presented some blow-up results for strong solutions in the case $\sigma=1, \mu=0$ and gave some global existence results to the system for $\sigma=1, \mu=0$.
Moreover, In \cite{cfg}, the authors showed that there are solitary waves with singularities, like peakons and
cuspons, depending on the values of the rotating parameter $\Omega$ and the balance index $\sigma$, they also
proved that horizontally symmetric weak solutions of this model must be traveling waves. In \cite{zl}, the authors established the local well-posedness in the critical Besov
space $B_{2,1}^{\frac{3}{2}}\times B_{2,1}^{\frac{1}{2}},$ presented a blow-up result with the initial data
in critical Besov space, studied the Gevrey regularity and analyticity of solutions to the system in a range of
Gevrey-Sobolev spaces in the sense of Hardamard and obtained a precise lower
bound of the lifespan.

The aim of this paper is to derive precise wave-breaking scenario for system (1.1)
in the case $\sigma>0$ , to give two blow-up results and corresponding blow-up rate result.
A notable difference in the blow-up analysis for (1.1)
stems from the cubic term $\Omega\rho(\rho u)_{x},$ which fails to be controlled by the conservation
laws. By following the method in \cite{c-l}, we get a uniform upper bound of $u_{x}$ even if we don¡¯t know whether it is bounded from below.

Our paper is organized as follows. In Section 2, some preliminary
estimates and results are recalled and presented. Section 3 is devoted to the proof
of precise wave-breaking scenario. In Section 4, we
provide a sufficient condition for blow-up solutions. Section 5 is about the blow-up rate of strong
solutions.  \\
\newline
\textbf{Notation}  Given a Banach space $Z$, we denote its norm by
 $\|\cdot\|_{Z}$. Since all space of functions are over
 $\mathbb{R}$, for simplicity, we drop $\mathbb{R}$ in our notations
  if there is no ambiguity.

\section{Preliminaries}
\newtheorem{theorem2}{Theorem}[section]
\newtheorem{lemma2}{Lemma}[section]
\newtheorem {remark2}{Remark}[section]
\newtheorem{corollary2}{Corollary}[section]

We now provide the framework in which we shall reformulate the
system (1.1). Note that if $p(x):=\frac{1}{2}e^{-|x|}$, $x\in
\mathbb{R}$, then $(1-
\partial^{2}_{x})^{-1}f = p\ast f $ for all $f \in L^{2}$. Here we denote by $\ast$ the convolution. Using
this two identities, we can rewrite the system (1.1) as follows:
\begin{align}\label{eq:original}
\left\{\begin{array}{ll}u_{t}+(\sigma u-\mu)u_{x}=-\partial_{x}p\ast((\mu-A)u+\frac{3-\sigma}{2}u^{2}+\frac{\sigma}{2}u_{x}^{2}+\frac{1-2\Omega A}{2}\rho^{2}\\
\ \ \ \ \ \ \ \ \ \ \ \ \ \ \ \ \ \ \ \ \ \ \ \ \ \ \ \ \ \ \  \ \ \ \ \ \ \ \ \ \ \  \ \ \ \ \ \ -\Omega\rho^{2}u)+\Omega p\ast(\rho^{2}u_{x}),\\
 \rho_{t}+u\rho_{x}=-\rho u_{x}, \\
u(0,x) =u_{0}(x), \\
\rho(0,x) = \rho_{0}(x).\end{array}\right.
\end{align}
\begin{lemma2}\cite{f-g-l}
Given $z_{0}=(u_{0},\rho_{0}-1)\in H^{s}\times
H^{s-1}, s> \frac{3}{2},$ there exists a maximal $T=T(\parallel
z_{0}\parallel_{H^{s}\times H^{s-1}})>0$ and a unique solution
$z=(u,\rho-1)$ to the system (2.1) such that
$$
z=z(\cdot,z_{0})\in C([0,T); H^{s}\times H^{s-1})\cap
C^{1}([0,T);H^{s-1}\times H^{s-2}).
$$
Moreover, the solution depends continuously on the initial data,
i.e. the mapping
$$z_{0}\rightarrow z(\cdot,z_{0}): H^{s}\times H^{s-1}\rightarrow
C([0,T); H^{s}\times H^{s-1})\cap C^{1}([0,T);H^{s-1}\times H^{s-2})
$$
is continuous.
\end{lemma2}
\begin{lemma2} Let $z_{0}=(u_{0},\rho_{0}-1)\in H^{s}\times
H^{s-1}, s> \frac{3}{2},$ and let $T>0$ be the
maximal existence time of the corresponding solution $z=(u,\rho-1)$
to system (2.1). Then we have
\begin{align*}E(t)&=\int_{\mathbb{R}} (u^2(t,x)+u_x^2(t,x)+(1-2\Omega A)(\rho-1)^{2}(t,x))dx\\
&=\int_{\mathbb{R}} (u_{0}^2(x)+u_{0,x}^2(x)+(1-2\Omega A)(\rho_{0}-1)^{2}(x))dx:=E(0).
\end{align*}
\end{lemma2}
\textbf{Proof}~ Multiplying the first equation of (1.1)by $u$ and integrating by parts, we have
\begin{equation}
\frac{1}{2}\frac{d}{dt}\int_{\mathbb{R}} (u^2(t,x)+u_x^2(t,x))dx=-(1-2\Omega A)\int_{\mathbb{R}}u \rho\rho_{x}dx.
\end{equation}
Rewrite the second equation in (1.1), we have
\begin{equation}
(\rho-1)_{t}+(u(\rho-1))_{x}+u_{x}=0.
\end{equation}
Using (2.3), a direct computation implies
\begin{align}&\frac{1}{2}\frac{d}{dt}\int_{\mathbb{R}}((1-2\Omega A))(\rho-1)^{2}dx\nonumber\\
=&(1-2\Omega A)\int_{\mathbb{R}}(\rho-1)(\rho-1)_{t}dx\nonumber\\
=&-(1-2\Omega A)\int_{\mathbb{R}}(\rho-1)((u(\rho-1))_{x}+u_{x})dx\nonumber\\
=&(1-2\Omega A)\int_{\mathbb{R}}u\rho \rho_{x}dx.
\end{align}
Adding equations (2.2) and (2.4), we have
$$\frac{d}{dt}\int_{\mathbb{R}} (u^2(t,x)+u_x^2(t,x)+(1-2\Omega A)(\rho-1)^{2})dx=0.$$
This completes the proof of the Lemma 2.2.
\begin{lemma2}\cite{f-g-l}Assume $(u,\rho)$ is the solution of system (2.1) with initial data $(u_{0}, \rho_{0}-1)\in H^{s}\times
H^{s-1} ,  s> \frac{3}{2},$  and let $T$ be the maximal time of existence. Then
$$T<+\infty \Rightarrow \int_{0}^{T}\|u_{x}(\tau)\|_{L^{\infty}}d\tau=+\infty.$$
\end{lemma2}
\begin{lemma2}\cite{f-g-l} Assume that $1-2\Omega A>0$. Let $(u_{0}, \rho_{0}-1)\in H^{s}\times
H^{s-1} \ \text{with} \ s> \frac{3}{2},$ and $T>0$ be the maximal time of existence of the solution $(u,\rho)$ to system (2.1)
with initial data $(u_{0}, \rho_{0}-1).$ Then the corresponding solution $(u,\rho)$ blows up in finite time $T<\infty$ if
and only if
$$\lim_{t\rightarrow T^{-}}\left\{\sup_{x\in \mathbb{R}}|u_{x}(t,x)|\right\}=+\infty.$$
\end{lemma2}
\ \ \ \ \ \ \ Consider now the following initial value problem
\begin{equation}
\left\{\begin{array}{ll}\frac{\partial q}{\partial t}=u(t, q),\ \ \ \ t\in[0,T), \\
q(0,x)=x,\ \ \ \ x\in\mathbb{R}, \end{array}\right.
\end{equation}
where $u\in C([0,T);H^s)\cap
C^1([0,T);H^{s-1})$ is the first component of the solution $z=(u,\rho)$ to (1.1).
Applying classical results in the theory of ordinary differential
equations, one can obtain two results on $q$ which are
crucial in studying  blow-up phenomena.
\begin{lemma2}\cite{f-g-l} Let $u\in C([0,T);H^s)\cap
C^1([0,T);H^{s-1}), s >\frac{3}{2}$. Then Eq.(2.5) has a unique solution
$q\in C^1([0,T)\times \mathbb{R};\mathbb{R})$. Moreover, the map
$q(t,\cdot)$ is an increasing diffeomorphism of $\mathbb{R}$ with
$$
q_{x}(t,x)=\exp\left(\int_{0}^{t}u_{x}(s,q(s,x))ds\right)>
0, \ \ \forall(t,x)\in [0,T)\times \mathbb{R}.$$
Consequently, the $L^{\infty}$-norm of any function $v(t, \cdot)\in L^{\infty}$ $(t\in [0, T))$ is preserved under the
family of the diffeomorphisms $q(t, \cdot),$  i.e.
$$\|v(t, \cdot)\|_{L^{\infty}}=\|v(t, q(t, \cdot))\|_{L^{\infty}}, \ \ \ \ t\in [0, T).$$
Similarly,
\begin{equation}
\inf_{x\in \mathbb{R}}v(t,x)= \inf_{x\in \mathbb{R}}v(t,q(t,x)), \ \ \ \ \ \sup_{x\in \mathbb{R}}v(t,x)= \sup_{x\in \mathbb{R}}v(t,q(t,x)),\ \ t\in [0, T).
\end{equation}
\end{lemma2}
\begin{lemma2}\cite{Constantin 4}
Let $T>0$ and $v\in C^{1}([0,T); H^{2})$. Then
for every $t\in[0,T)$ there exists at least one point $\xi(t)\in
\mathbb{R}$ with
$$ m(t):=\inf_{x\in \mathbb{R}}[v_{x}(t,x)]=v_{x}(t,\xi(t)),$$ and
the function $m$ is almost everywhere differentiable on $(0,t_{0})$
with $$ \frac{d}{dt}m(t)=v_{tx}(t,\xi(t)) \ \ \ \ a.e.\ \ on \
(0,T).$$
\end{lemma2}
\section{Precise blow-up scenarios}
\newtheorem{theorem3}{Theorem}[section]
\newtheorem{lemma3}{Lemma}[section]
\newtheorem {remark3}{Remark}[section]
\newtheorem{corollary3}{Corollary}[section]

In this section, we will derive the precise blow-up scenario for
strong solutions to the system (1.1) with $\sigma>0$, which improves the corresponding result for $\sigma=1$ and $\mu=0$ in \cite{f-g-l}. Moreover, in
\cite{cfg}, the authors proved that $u_{x}$ is uniformly bounded from above on the set $[0,T)\times \Lambda,$ where $\Lambda=\{x\in \mathbb{R}: \rho_{0}(x)=0\},$
even if one don't know whether it is bounded from below. Here, we improve the estimate to all of $\mathbb{R},$ which is useful to the prove of the precise blow-up scenario.

\begin{lemma3} Suppose that $1-2\Omega A>0$ and $\sigma>0.$ Let $(u, \rho)$ be the solution of (2.1) with initial data $(u_{0},\rho_{0}-1)\in H^{s}\times
H^{s-1}, s> \frac{3}{2},$ and $T$ be the maximal time of existence. Then
\begin{equation}
\sup_{x\in \mathbb{R}}u_{x}(t,x)\leq \|u_{0,x}\|_{L^{\infty}}+\sqrt{\frac{(1-2\Omega A)\|\rho_{0}\|_{L^{\infty}}^{2}+C^{2}}{\sigma}},
\end{equation}
where $C$ is a positive constant depend on $E(0)$  and $\|\rho_{0}\|_{L^{\infty}},$ it will be given in (3.13) below.
\end{lemma3}
\textbf{Proof}~ The local well-posedness theorem and a density argument imply
that it suffices to prove the desired estimates for $s\geq 3$. Thus, we take $s=3$ in the proof.
Differentiating the first equation in (2.1) with respect to $x$ and using the identity $-\partial^{2}_{x}p\ast f = f- p\ast f,$
we have
\begin{equation}
u_{tx}+(\sigma u-\mu)u_{xx}=-\frac{\sigma}{2}u_{x}^{2}+\frac{1-2\Omega A}{2}\rho^{2}+f(t,x),
\end{equation}
where
$$f(t,x)=-(\mu-A)\partial_{x}^{2}p\ast u+\frac{3-\sigma}{2}u^{2}-\Omega\rho^{2}u-p\ast(\frac{3-\sigma}{2}u^{2}+\frac{\sigma}{2}u_{x}^{2}+\frac{1-2\Omega A}{2}\rho^{2}-\Omega\rho^{2}u)+
\Omega\partial_{x}p\ast(\rho^{2}u_{x}).$$
Let $$M(t)=\sup_{x\in\mathbb{R}}u_{x}(t,x),\ \ \ \ \ t\in[0,T).$$
By Lemma 2.6 and the fact $$\sup_{x\in\mathbb{R}}u_{x}(t,x)=-\inf_{x\in\mathbb{R}}(-u_{x}(t,x)),$$
there exists a point $\xi(t)\in \mathbb{R}$ such that $M(t)=\sup\limits_{x\in\mathbb{R}}u_{x}(t,x)=u_{x}(t,\xi(t)).$ Obviously, $u_{xx}(t,\xi(t))=0$ and $M^{\prime}(t)=u_{tx}(t,\xi(t)).$
Take the trajectory $q(t,x)$ defined in (2.5). By Lemma 2.5 we have $q(t,\cdot)$ is an increasing diffeomorphism of $\mathbb{R}$
for every $t\in [0,T).$ Therefore, there exists $x_{1}(t)\in \mathbb{R}$ such that
$\xi(t)=q(t, x_{1}(t)).$  Let $$\gamma(t)=\rho(t, q(t, x_{1}(t)), \ \ \ \ t\in[0,T).$$
By (3.2), the second equation in (2.1) and (2.5), we can obtain
\begin{equation}
M^{\prime}(t)=-\frac{\sigma}{2}M^{2}(t)+\frac{1-2\Omega A}{2}\gamma^{2}(t)+f(t,q(t,x_{1}(t)))
\end{equation}
and
\begin{equation}
\gamma^{\prime}(t)=-M(t)\gamma(t),
\end{equation}
for $t\in[0,T),$  where $\prime$ denotes the derivative with respect to $t$. Note that $\partial_{x}^{2}p\ast u=\partial_{x}p\ast \partial_{x} u,$ we can rewrite $f(t,q(t,x_{1}(t)))$
as follows
\begin{align}
&\left|f(t,q(t,x_{1}(t)))\right|\nonumber\\
=&|(A-\mu)\partial_{x}p\ast \partial_{x} u+\frac{3-\sigma}{2}u^{2}-\Omega\rho^{2}u-p\ast(\frac{3-\sigma}{2}u^{2}+\frac{\sigma}{2}u_{x}^{2})-\frac{1-2\Omega A}{2}p\ast(\rho-1)^{2}\nonumber\\
&-(1-2\Omega A)p\ast(\rho-1)-\frac{1-2\Omega A}{2}+\Omega p\ast((\rho-1)^{2}u)+2 \Omega p\ast(\rho-1)u+\Omega p\ast u\nonumber\\
&+\Omega \partial_{x}p\ast(\rho(\rho-1)u_{x})+\Omega\partial_{x}p\ast((\rho-1)u_{x})+\Omega\partial_{x}p\ast u_{x}|\nonumber\\
\leq&|A-\mu| |\partial_{x}p\ast \partial_{x} u|+|\frac{3-\sigma}{2}u^{2}|+\Omega|\rho^{2}u|+|p\ast(\frac{3-\sigma}{2}u^{2}+\frac{\sigma}{2}u_{x}^{2})|+\frac{1-2\Omega A}{2}|p\ast(\rho-1)^{2}|\nonumber\\
&+(1-2\Omega A)|p\ast(\rho-1)|+\frac{1-2\Omega A}{2}+\Omega |p\ast((\rho-1)^{2}u)|+2 \Omega |p\ast(\rho-1)u|+\Omega| p\ast u|\nonumber\\
&+\Omega |\partial_{x}p\ast[\rho(\rho-1)u_{x}]|+\Omega|\partial_{x}p\ast((\rho-1)u_{x})|+\Omega|\partial_{x}p\ast u_{x}|
\end{align}
Since now $s=3,$ we have $u\in C_{0}^{1},$ subscript $0$ means the function decays to zero at infinity. It follows that
$$\inf_{x\in\mathbb{R}}u_{x}(t,x)\leq 0 , \ \ \ \ \sup_{x\in\mathbb{R}}u_{x}(t,x)\geq 0,\ \ \ t\in[0,T). $$
Then, $M(t)\geq 0.$ From equation (3.4), a direct computation implies
$$|\gamma (t)|=|\gamma (0)|e^{\int_{0}^{t}-M(\tau)d\tau},\ \ \ t\in[0,T).$$
So
\begin{equation}
|\rho(t, q(t, x_{1}(t))|=|\gamma (t)|\leq|\gamma (0)|=|\rho_{0}(x_{1}(0))|\leq \|\rho_{0}\|_{L^{\infty}},\ \ \ t\in[0,T).
\end{equation}
Next, we will estimate (3.5) item by item, here we will use (3.6), $Young$ inequality, $H\ddot{o}lder$ inequality and the fact $\|\partial_{x}p\|_{L^{2}}=\frac{1}{2}=\|p\|_{L^{\infty}}.$
$$|A-\mu| |\partial_{x}p\ast \partial_{x} u|\leq |A-\mu| \|\partial_{x}p\|_{L^{2}}\|u_{x}\|_{L^{2}}=\frac{|A-\mu|}{2}\|u_{x}\|_{L^{2}}\leq\frac{1}{4}+\frac{(A-\mu)^{2}}{4}\|u_{x}\|_{L^{2}}^{2},$$
$$|\frac{3-\sigma}{2}u^{2}|\leq \frac{|3-\sigma|}{2}\|u\|_{L^{\infty}}^{2}\leq\frac{|3-\sigma|}{4}\|u\|_{H^{1}}^{2}=\frac{|3-\sigma|}{4}(\|u\|_{L^{2}}^{2}+\|u_{x}\|_{L^{2}}^{2}),$$
$$\Omega|\rho^{2}u|\leq \Omega|\rho|^{2}\|u\|_{L^{\infty}}\leq\frac{\Omega^{2}}{2}|\rho|^{4}+\frac{1}{4}(\|u\|_{L^{2}}^{2}+\|u_{x}\|_{L^{2}}^{2})\leq \frac{\Omega^{2}}{2}\|\rho_{0}\|_{L^{\infty}}^{4}+\frac{1}{4}(\|u\|_{L^{2}}^{2}+\|u_{x}\|_{L^{2}}^{2}),$$
$$\left|p\ast(\frac{3-\sigma}{2}u^{2}+\frac{\sigma}{2}u_{x}^{2})\right|\leq \|p\|_{L^{\infty}}\left\|\frac{3-\sigma}{2}u^{2}+\frac{\sigma}{2}u_{x}^{2}\right\|_{L^{1}}=\frac{|3-\sigma|}{4}\|u\|_{L^{2}}^{2}+\frac{|\sigma|}{4}\|u_{x}\|_{L^{2}}^{2},$$
$$\frac{1-2\Omega A}{2}|p\ast(\rho-1)^{2}|\leq \frac{1-2\Omega A}{2}\|p\ast(\rho-1)^{2}\|_{L^{\infty}}\leq \frac{1-2\Omega A}{2}\|p\|_{L^{\infty}}\|(\rho-1)^{2}\|_{L^{1}}=\frac{1-2\Omega A}{4}\|\rho-1\|_{L^{2}}^{2},$$
$$(1-2\Omega A)|p\ast(\rho-1)|\leq (1-2\Omega A)\|p\|_{L^{2}}\|\rho-1\|_{L^{2}}\leq \frac{1-2\Omega A}{4}+\frac{1-2\Omega A}{4}\|\rho-1\|_{L^{2}}^{2},$$
\begin{equation}
\Omega |p\ast((\rho-1)^{2}u)|\leq \frac{\Omega}{2}\|(\rho-1)^{2}u\|_{L^{1}}\leq \frac{\Omega}{2}\|u\|_{L^{\infty}}\|\rho-1\|_{L^{2}}^{2},
\end{equation}
\begin{equation}
2\Omega|p\ast((\rho-1)u)|\leq2\Omega\|p\|_{L^{2}}\|(\rho-1)u\|_{L^{2}}\leq\Omega\|(\rho-1)u\|_{L^{2}}\leq\Omega\|u\|_{L^{\infty}}\|\rho-1\|_{L^{2}},
\end{equation}
\begin{equation}
\Omega|p\ast u|\leq\Omega\|p\|_{L^{2}}\|u\|_{L^{2}}\leq\frac{1}{4}+\frac{\Omega^{2}}{4}\|u\|_{L^{2}}^{2},
\end{equation}
\begin{align}
\Omega |\partial_{x}p\ast[\rho(\rho-1)u_{x}]|\leq \Omega\|\partial_{x}p\|_{L^{\infty}}\|\rho(\rho-1)u_{x}\|_{L^{1}}\leq\frac{\Omega}{2}\|\rho_{0}\|_{L^{\infty}}\|\rho-1\|_{L^{2}}\|u_{x}\|_{L^{2}}\nonumber\\
\leq\frac{\Omega}{4}\|\rho_{0}\|_{L^{\infty}}(\|\rho-1\|_{L^{2}}^{2}+\|u_{x}\|_{L^{2}}^{2}),
\end{align}
\begin{equation}
\Omega |\partial_{x}p\ast((\rho-1)u_{x})|\leq\frac{\Omega}{2}\|(\rho-1)u_{x}\|_{L^{1}}\leq\frac{\Omega}{4}(\|\rho-1\|_{L^{2}}^{2}+\|u_{x}\|_{L^{2}}^{2}),
\end{equation}
\begin{equation}
\Omega|\partial_{x}p\ast u_{x}|\leq\frac{\Omega}{2}\| u_{x}\|_{L^{2}}\leq\frac{1}{4}+\frac{\Omega^{2}}{4}\|u_{x}\|_{L^{2}}^{2}.
\end{equation}
Since
$$\|u\|_{L^{2}}^{2}\leq E(0),$$ $$\|u_{x}\|_{L^{2}}^{2}\leq E(0),$$ $$\|\rho-1\|_{L^{2}}^{2}\leq \frac{1}{1-2\Omega A}E(0),$$
$$\|u\|_{L^{\infty}}^{2}\leq \frac{1}{2}\|u\|_{H^{1}}^{2}\leq \frac{1}{2}E(0),$$
we have
\begin{align}
|f|\leq \frac{3(1-\Omega A)}{2}+\left[\frac{(A-\mu)^{2}}{4}+\frac{|3-\sigma|}{2}+\frac{3}{4}+\frac{|\sigma|}{4}+\frac{\Omega}{2}+\frac{3\Omega}{4(1-2\Omega A)}+\frac{\Omega^{2}}{2}\right]E(0)\nonumber\\
+\frac{\Omega^{2}}{2}\|\rho_{0}\|_{L^{\infty}}^{4}+\left(\frac{\Omega}{4(1-2\Omega A)}+\frac{\Omega}{4}\right)\|\rho_{0}\|_{L^{\infty}}E(0)+\frac{\sqrt{2}\Omega}{4(1-2\Omega A)}\sqrt{E(0)}E(0):=\frac{C^{2}}{2}
\end{align}
where constant $C=C(E(0),\ \|\rho_{0}\|_{L^{\infty}})>0.$
For any given $x\in \mathbb{R},$ define
$$ P(t)= M(t)-\|u_{0,x}\|_{L^{\infty}}-\sqrt{\frac{(1-2\Omega A)\|\rho_{0}\|_{L^{\infty}}^{2}+C^{2}}{\sigma}}.$$
Observe $P(t)$ is a $C^{1}$-differentiable function in $[0,t)$ and satisfies
\begin{align*}
P(0)= M(0)-\|u_{0,x}\|_{L^{\infty}}-\sqrt{\frac{(1-2\Omega A)\|\rho_{0}\|_{L^{\infty}}^{2}+C^{2}}{\sigma}}\leq  M(0)-\|u_{0,x}\|_{L^{\infty}}\\
=\sup_{x\in \mathbb{R}}u_{0,x}(x)-\|u_{0,x}\|_{L^{\infty}}\leq 0.
\end{align*}
Next, we will prove that $P(t)\leq 0$ for $t\in[0,T).$ If not, then there is a $t_{0}\in [0,T)$ such that $P(t_{0})>0.$ Let
$$t_{1}=\max\{t<t_{0};P(t)=0\}.$$
It follows that $P(t_{1})=0$ and $P^{\prime}(t_{1})\geq0,$ or equivalently,
$$M(t_{1})=\|u_{0,x}\|_{L^{\infty}}+\sqrt{\frac{(1-2\Omega A)\|\rho_{0}\|_{L^{\infty}}^{2}+C^{2}}{\sigma}}$$
and $M^{\prime}(t_{1})\geq0.$
On the other hand, by (3.3), (3.6) and (3.13) we get
\begin{align*}
M^{\prime}(t_{1})&=-\frac{\sigma}{2}M^{2}(t_{1})+\frac{1-2\Omega A}{2}\gamma^{2}(t_{1})+f(t,q(t,x_{1}(t_{1})))\\
&\leq -\frac{\sigma}{2}\left(\|u_{0,x}\|_{L^{\infty}}+\sqrt{\frac{(1-2\Omega A)\|\rho_{0}\|_{L^{\infty}}^{2}+C^{2}}{\sigma}}\right)^{2}+\frac{1-2\Omega A}{2}\|\rho_{0}\|_{L^{\infty}}^{2}+\frac{C^{2}}{2}\\
&< -\frac{\sigma}{2}\|u_{0,x}\|_{L^{\infty}}^{2}-\frac{\sigma}{2}\cdot \frac{(1-2\Omega A)\|\rho_{0}\|_{L^{\infty}}^{2}+C^{2}}{\sigma}+\frac{1-2\Omega A}{2}\|\rho_{0}\|_{L^{\infty}}^{2}+\frac{C^{2}}{2}\\
&= -\frac{\sigma}{2}\|u_{0,x}\|_{L^{\infty}}^{2}\leq 0,
\end{align*}
a contradiction. Therefore,  $P(t)\leq 0$ for $t\in[0,T).$ Since $x$ is chosen arbitrarily, we obtain (3.1).  This completes the proof of Lemma 3.1.

Our next result describes the precise blow-up scenario for
sufficiently regular solutions to system (2.1).
\begin{theorem3}
Assume that $1-2\Omega A>0$ and $\sigma>0$. Let $(u_{0}, \rho_{0}-1)\in H^{s}\times
H^{s-1} \ \text{with} \ s> \frac{3}{2},$ and $T>0$ be the maximal time of existence of the solution $(u,\rho)$ to system (2.1)
with initial data $(u_{0}, \rho_{0}-1).$ Then the corresponding solution $(u,\rho)$ blows up in finite time $T<+\infty$ if
and only if
$$\lim_{t\rightarrow T}\left\{\inf_{x\in \mathbb{R}}u_{x}(t,x)\right\}=-\infty.$$
\end{theorem3}
\textbf{Proof}~ By Lemma 2.1  and a density argument imply
that it suffices to prove the desired estimates for $s\geq 3$. Thus, we take $s=3$ in the proof.

By Sobolev's imbedding theorem $H^{s}\hookrightarrow L^{\infty}$ with $s>\frac{1}{2},$ it is clear that if
$$\lim_{t\rightarrow T}\left\{\inf_{x\in \mathbb{R}}u_{x}(t,x)\right\}=-\infty,$$ then $T<+\infty$.

Conversely, let $T<+\infty$ and assume that for some constant $K>0$ such that
$$u_{x}(t,x)\geq -K, \ \ \ \ \forall(t,x)\in [0,T)\times \mathbb{R}.$$ Then, it follows from Lemma 3.1 that
$|u_{x}(t,x)|\leq K_{1},$ where $K_{1}=K_{1}(K, E(0),\ \|\rho_{0}\|_{L^{\infty}}).$ Therefore, Lemma 2.3 ensures that the
maximal existence time $T=+\infty$, which contradicts the assumption that $T<+\infty.$ This completes the proof of Theorem 3.1.\\

\section{Blow-up}
\newtheorem{theorem4}{Theorem}[section]
\newtheorem{lemma4}{Lemma}[section]
\newtheorem {remark4}{Remark}[section]
\newtheorem{corollary4}{Corollary}[section]

In this section, we discuss the blow-up phenomena of system (2.1) and
prove that there exist strong solutions to system (2.1) which do not
exist globally in time. Note that estimate (3.13) is also true for $\sigma<0$, we will use (3.13) in this section and the following section.
\begin{theorem4}Assume that $1-2\Omega A>0$ and $\sigma<0$. Let $(u_{0}, \rho_{0}-1)\in H^{s}\times
H^{s-1} \ \text{with} \ s> \frac{3}{2}$ and $T>0$ be the maximal time of existence of the solution $(u,\rho)$ to system (2.1)
with initial data $(u_{0}, \rho_{0}-1).$  If there exists some $x_{0}\in \mathbb{R}$ such that $u_{0}^{\prime}(x_{0})>\frac{C}{\sqrt{-\sigma}},$ where $C>0$ is defined in
(3.13). Then the corresponding solution to system (2.1)
blows up in finite time in the following sense: there exists a $T_{1}$ with
$$0<T_{1}\leq\frac{-2}{\sigma u_{0}^{\prime}(x_{0})-\sqrt{C(-\sigma)^{\frac{3}{2}}u_{0}^{\prime}(x_{0})}}$$ such that
\begin{equation}
\liminf_{t\rightarrow T_{1}}\sup_{x\in \mathbb{R}}u_{x}(t,x)=+\infty.
\end{equation}
\end{theorem4}
\textbf{Proof}~ Similar to the proof of Lemma 3.1, we consider the functions
$$M(t)=\sup\limits_{x\in\mathbb{R}}u_{x}(t,x)=u_{x}(t,\xi(t))$$
and
$$\gamma(t)=\rho(t, q(t, x_{1}(t)), $$
where $\xi(t)=q(t, x_{1}(t)).$
By (3.3) and (3.13), we have
\begin{equation}
M^{\prime}(t)=-\frac{\sigma}{2}M^{2}(t)+\frac{1-2\Omega A}{2}\gamma^{2}(t)+f(t,q(t,x_{1}(t)))\geq-\frac{\sigma}{2}M^{2}(t)-\frac{C^{2}}{2}
\end{equation}
Since there exists some $x_{0}\in \mathbb{R}$ such that $u_{0}^{\prime}(x_{0})>\frac{C}{\sqrt{-\sigma}},$
so $$ M(0)=\sup_{x\in \mathbb{R}}u_{0}^{\prime}(x)\geq u_{0}^{\prime}(x_{0})>\frac{C}{\sqrt{-\sigma}},$$ it follows that
$M^{\prime}(0)\geq-\frac{\sigma}{2}M^{2}(0)-\frac{C^{2}}{2}>0.$ Then, we will prove that
\begin{equation}
M(t)>\frac{C}{\sqrt{-\sigma}}  \ \ \ \text{for} \ \ \ t\in[0,T).
\end{equation}
If not, then suppose there is a $t_{0}\in [0,T)$ such that $M(t_{0})\leq\frac{C}{\sqrt{-\sigma}}.$ Define
$$t_{1}=\max\{t<t_{0};M(t)=\frac{C}{\sqrt{-\sigma}}\}.$$
It follows that $M(t_{1})=\frac{C}{\sqrt{-\sigma}}$ and $M^{\prime}(t_{1})< 0.$ On the other hand, by (4.2) we have
$$M^{\prime}(t_{1})\geq-\frac{\sigma}{2}M^{2}(t_{1})-\frac{C^{2}}{2}=0.$$
This contradiction implies (4.3). Using (4.2) and (4.3), we see that $M(t)$  is strictly
increasing over $[0, T).$  Therefore $M(t)>M(0)\geq u_{0}^{\prime}(x_{0})>\frac{C}{\sqrt{-\sigma}}.$ Let
 $$ \delta=\frac{1}{2}+\frac{1}{2}\sqrt{\frac{C}{u_{0}^{\prime}(x_{0})\sqrt{-\sigma}}}\in \left(\frac{1}{2},1\right).$$
By (4.2), we know
\begin{align*}
M^{\prime}(t)\geq-\frac{\sigma}{2}M^{2}(t)-\frac{C^{2}}{2}&=-\frac{\sigma}{2}M^{2}(t)\left[1+\frac{C^{2}}{\sigma M^{2}(t)}\right]\\
&\geq-\frac{\sigma}{2}M^{2}(t)[1-(2\delta-1)^{4}]\\
&\geq -\frac{\sigma}{2}M^{2}(t)\cdot 2\delta=-\sigma\delta M^{2}(t).
\end{align*}
Solving this inequality we obtain
$$\sup_{x\in\mathbb{R}}u_{x}(t,x)=M(t)\geq \frac{u_{0}^{\prime}(x_{0})}{1+\sigma\delta u_{0}^{\prime}(x_{0})t}\rightarrow +\infty \ \ \ \ \text{as} \ \ t\rightarrow -\frac{1}{\sigma\delta u_{0}^{\prime}(x_{0})}=\frac{-2}{\sigma u_{0}^{\prime}(x_{0})-\sqrt{C(-\sigma)^{\frac{3}{2}}u_{0}^{\prime}(x_{0})}},$$
which proves (4.1). By Lemma 2.4, we have the corresponding solution $(u,\rho)$ blows up in finite time.\\

By Sobolev's imbedding theorem $H^{s}\hookrightarrow L^{\infty}$ with $s>\frac{1}{2},$ it is clear that if $\rho(t,x)$ becomes unbounded in finite time, then
the maximal time of existence of the solution $(u,\rho)$ to system (2.1) $T<+\infty,$ it follows from Theorem 3.1 that $u_{x}(t,x)$
must be unbounded from below in finite time. In order to study the fine structure of finite time singularities
we shall assume in the following that there is a $M>0$ such that $\|\rho(t,\cdot)\|_{L^{\infty}}\leq M$ for all $t\in [0,T).$ Without loss of generality, we let
$\sigma=1,$ $\mu=0$ in the following result.
\begin{theorem4}
Assume that $1-2\Omega A>0$ and $\sigma=1,$ $\mu=0$. Let $(u_{0}, \rho_{0}-1)\in H^{s}\times
H^{s-1} \ \text{with} \ s> \frac{3}{2}$ and $T>0$ be the maximal time of existence of the solution $(u,\rho)$ to system (2.1)
with initial data $(u_{0}, \rho_{0}-1).$ Assume
that there exists $M>0$ such that $\|\rho(t,\cdot)\|_{L^{\infty}}\leq M$ for all $t\in [0,T).$ If
\begin{equation}
\int_{\mathbb{R}}u_{0,x}^{3}dx\leq -\sqrt{2E(0)N},
\end{equation}
where
\begin{align}
N=&\left(\frac{3M^{2}(1-2\Omega A)}{2}+\frac{9}{4}\right)E(0)+\frac{3\sqrt{2}\Omega M^{2}}{2}(E(0))^{\frac{3}{2}}+\frac{3\sqrt{2}\Omega}{4(1-2\Omega A)}(E(0))^{\frac{5}{2}}\nonumber\\
&+\left[\frac{6+3A^{2}+6\Omega^{2}}{4}+\frac{3\sqrt{2}\Omega}{2\sqrt{1-2\Omega A}}+\frac{3\Omega(1-\Omega A)(M+1)}{2(1-2\Omega A)}\right]E(0)^{2}\nonumber,
\end{align}
then the corresponding solution to system (2.1) blows up
in finite time.
\end{theorem4}
\textbf{Proof}  As mentioned earlier, here we only need
to show that the above theorem holds for $s=3$. Let $\sigma=1$ and $\mu=0$ in (3.2), we have
\begin{align*}
u_{tx}+uu_{xx}=&-\frac{1}{2}u_{x}^{2}+\frac{1-2\Omega A}{2}\rho^{2}+u^{2}-\Omega u\rho^{2}+A\partial_{x}^{2}p\ast u\\
&-p\ast(u^{2}+\frac{1}{2}u_{x}^{2}+\frac{1-2\Omega A}{2}\rho^{2}-\Omega\rho^{2}u)+\Omega\partial_{x}p\ast(\rho^{2}u_{x}).
\end{align*}
Then, it follws that
\begin{align*}
&\frac{d}{dt}\int_{\mathbb{R}}u_{x}^{3}dx\\
=&\int_{\mathbb{R}}3u_{x}^{2}u_{tx}dx\\
=&3\int_{\mathbb{R}}u_{x}^{2}(-uu_{xx}-\frac{1}{2}u_{x}^{2}+\frac{1-2\Omega A}{2}\rho^{2}+u^{2}-\Omega u\rho^{2}+A\partial_{x}^{2}p\ast u\\
&-p\ast(u^{2}+\frac{1}{2}u_{x}^{2}+\frac{1-2\Omega A}{2}\rho^{2}-\Omega\rho^{2}u)+\Omega\partial_{x}p\ast(\rho^{2}u_{x}))dx\\
\leq &-\frac{1}{2}\int_{\mathbb{R}}u_{x}^{4}dx+3\int_{\mathbb{R}}u_{x}^{2}(\frac{1-2\Omega A}{2}\rho^{2}+u^{2}-\Omega u\rho^{2}+A\partial_{x}^{2}p\ast u+\Omega p\ast(\rho^{2}u)+\Omega\partial_{x}p\ast(\rho^{2}u_{x}))dx,
\end{align*}
here we used $-3\int_{\mathbb{R}}uu_{x}^{2}u_{xx}dx=\int_{\mathbb{R}}u_{x}^{4}dx.$
Since
$$\|u\|_{L^{2}}^{2}\leq E(0),$$ $$\|u_{x}\|_{L^{2}}^{2}\leq E(0),$$ $$\|\rho-1\|_{L^{2}}^{2}\leq \frac{1}{1-2\Omega A}E(0),$$
$$\|u\|_{L^{\infty}}^{2}\leq \frac{1}{2}\|u\|_{H^{1}}^{2}\leq \frac{1}{2}E(0),$$
we get
$$\frac{3(1-2\Omega A)}{2}\int_{\mathbb{R}}u_{x}^{2}\rho^{2}dx\leq\frac{3M^{2}(1-2\Omega A)}{2}E(0),$$
$$3\int_{\mathbb{R}}u_{x}^{2}u^{2}dx\leq 3\|u\|_{L^{\infty}}^{2}\int_{\mathbb{R}}u_{x}^{2}dx\leq \frac{3}{2}(E(0))^{2},$$
$$-3\Omega \int_{\mathbb{R}}u_{x}^{2}u\rho^{2}dx\leq3\Omega\|u\|_{L^{\infty}}\|\rho\|_{L^{\infty}}^{2}\int_{\mathbb{R}}u_{x}^{2}dx\leq \frac{3\sqrt{2}}{2}\Omega M^{2}\sqrt{E(0)}E(0),$$
$$3A\int_{\mathbb{R}}u_{x}^{2}\partial_{x}^{2}p\ast udx\leq 3|A|\|\partial_{x}p\ast \partial_{x} u\|_{L^{\infty}}\int_{\mathbb{R}}u_{x}^{2}dx\leq \frac{3}{4}E(0)+\frac{3A^{2}}{4}(E(0))^{2}.$$
Note that
$$\Omega p\ast\rho^{2}u=\Omega p\ast((\rho-1)^{2}u)+2\Omega p\ast((\rho-1)u)+\Omega p\ast u$$
and
$$\Omega\partial_{x}p\ast(\rho^{2}u_{x})=\Omega\partial_{x}p\ast[\rho(\rho-1)u_{x}]+\Omega\partial_{x}p\ast((\rho-1)u_{x})+\Omega \partial_{x}p\ast u_{x},$$
a similar estimation as in (3.7)-(3.12) yields
$$3\Omega \int_{\mathbb{R}}u_{x}^{2}p\ast\rho^{2}udx\leq \frac{3\sqrt{2}\Omega}{4(1-2\Omega A)}(E(0))^{\frac{5}{2}}+\frac{3\sqrt{2}\Omega}{2\sqrt{1-2\Omega A}}(E(0))^{2}+\frac{3}{4}E(0)+\frac{3\Omega^{2}}{4}(E(0))^{2},$$
$$3\Omega \int_{\mathbb{R}}u_{x}^{2}\partial_{x}p\ast(\rho^{2}u_{x})dx\leq \frac{3\Omega (1-\Omega A)(M+1)}{2(1-2\Omega A)}(E(0))^{2}+\frac{3}{4}E(0)+\frac{3\Omega^{2}}{4}(E(0))^{2},$$
Thus
\begin{align*}
\frac{d}{dt}\int_{\mathbb{R}}u_{x}^{3}dx\leq-\frac{1}{2}\int_{\mathbb{S}}u_{x}^{4}dx+N.
\end{align*}
Using the following inequality
$$
\left|\int_{\mathbb{R}}u_{x}^{3}dx\right| \leq
\left(\int_{\mathbb{R}}u_{x}^{4}dx\right)^{\frac{1}{2}}\left(\int_{\mathbb{R}}u_{x}^{2}dx\right)^{\frac{1}{2}}\leq\left(\int_{\mathbb{R}}u_{x}^{4}dx\right)^{\frac{1}{2}}
\sqrt{E(0)},$$
and letting
$$m(t)=\int_{\mathbb{R}}u_{x}^{3}dx,$$
 we have
\begin{align*}
\frac{d}{dt}m(t)\leq&-\frac{1}{2E(0)}m^{2}(t)+N\\
=&-\frac{1}{2E(0)}\left(m(t)+\sqrt{2E(0)N}\right)\left(m(t)-\sqrt{2E(0)N}\right)
\end{align*}
Note that if $m(0)<-\sqrt{2E(0)N}$ then
$m(t)<-\sqrt{2E(0)N}$ for all
$t\in[0,T)$. From the above inequality we obtain
\begin{eqnarray}
\nonumber \frac{m(0)+\sqrt{2E(0)N}}{m(0)-\sqrt{2E(0)N}}e^{\sqrt{\frac{2N}{E(0)}}t}-1\leq
\frac{2\sqrt{2E(0)N}}{m(t)-\sqrt{2E(0)N}}\leq0.
\end{eqnarray}
Since $0<\frac{m(0)+\sqrt{2E(0)N}}{m(0)-\sqrt{2E(0)N}}<1$,
then there exists
$$0<T\leq\sqrt{\frac{E(0)}{2N}}\ln\frac{m(0)-\sqrt{2E(0)N}}{m(0)+\sqrt{2E(0)N}},$$
such that $\lim_{t\rightarrow T} m(t)=-\infty.$ On the other hand,
$$|\int_{\mathbb{R}}u_{x}^{3}dx|\leq\|u_{x}\|_{L^{\infty}}\int_{\mathbb{R}}u_{x}^{2}dx\leq \|u_{x}\|_{L^{\infty}}\|u\|_{H^{1}}^{2}
\leq \|u_{x}\|_{L^{\infty}}E(0).$$
Applying Lemma 2.4, the corresponding solution to system (2.1) blows up
in finite time.

\section{Blow-up rate}
\newtheorem{theorem5}{Theorem}[section]
\newtheorem{lemma5}{Lemma}[section]
\newtheorem {remark5}{Remark}[section]
\newtheorem{corollary5}{Corollary}[section]
We now give more insight into the blow-up
mechanism for the wave-breaking solution to system (2.1).
\begin{theorem5}
Assume that $1-2\Omega A>0$ and $\sigma<0$. Let $(u_{0}, \rho_{0}-1)\in H^{s}\times
H^{s-1} \ \text{with} \ s> \frac{3}{2}$ and $T>0$ be the maximal time of existence of the solution $(u,\rho)$ to system (2.1)
with initial data $(u_{0}, \rho_{0}-1).$  If $T$ is finite in the following sense:
\begin{equation}
\lim_{t\rightarrow T}\sup_{x\in \mathbb{R}}u_{x}(t,x)=+\infty.
\end{equation}
Then
$$\lim_{t\rightarrow T}[(T-t)\sup_{x\in \mathbb{R}}u_{x}(t,x)]=-\frac{2}{\sigma}.$$
\end{theorem5}
\textbf{Proof} Similar to the proof of Lemma 3.1, we consider the functions
$$M(t)=\sup\limits_{x\in\mathbb{R}}u_{x}(t,x)=u_{x}(t,\xi(t))$$
and
$$\gamma(t)=\rho(t, q(t, x_{1}(t)), $$
where $\xi(t)=q(t, x_{1}(t)).$
By (3.3), we have
\begin{equation}
M^{\prime}(t)+\frac{\sigma}{2}M^{2}(t)=\frac{1-2\Omega A}{2}\gamma^{2}(t)+f(t,q(t,x_{1}(t)))
\end{equation}
Using (3.6) and (3.13), we have
$$\left|\frac{1-2\Omega A}{2}\gamma^{2}(t)+f(t,q(t,x_{1}(t)))\right|\leq \frac{1-2\Omega A}{2}\|\rho_{0}\|_{L^{\infty}}^{2}+\frac{C^{2}}{2}:=K_{2} \ \ \ a.e.\ \ on\ (0,T).$$
It follows that
\begin{equation}
-K_{2}\leq\frac{dM(t)}{dt}+\frac{\sigma}{2}M^{2}(t)\leq K_{2}.
\end{equation}
Let $\varepsilon\in(0,-\frac{\sigma}{2})$. Since $\lim\limits_{t\rightarrow
T}M(t)= +\infty$ by (5.1), there is some $t_{2}\in (0,T)$ with
$M^{2}(t_{2})>\frac{K_{2}}{\varepsilon}$. Let us first
prove that
\begin{equation}
M^{2}(t)>\frac{K_{2}}{\varepsilon}, \ \ \ \ t\in [t_{2}, T).
\end{equation}

Since $M$ is locally Lipschitz (it belongs to
$W_{loc}^{1,\infty}(\mathbb{R})$ by Lemma 2.6) there is some
$\theta>0$ such that
$$M^{2}(t)>\frac{K_{2}}{\varepsilon}, \ \ \ \ t\in [t_{2},
t_{2}+\theta).$$ Pick $\theta>0$ maximal with this property. If
$\theta < T-t_{2}$ we would have
$M^{2}(t_{2}+\theta)=\frac{K_{2}}{\varepsilon}$ while$$
\frac{dM}{dt}\geq-\frac{\sigma}{2}M^{2}(t)-K_{2}>-\frac{\sigma}{2}M^{2}(t)-\varepsilon
M^{2}(t)>0 \ \ \ \ \ \ a.e.\ on\ (t_{2}, t_{2}+\theta).$$ Being locally
Lipschitz, the function $M$ is absolutely continuous and therefore
we would obtain by integrating the previous relation on $[t_{2},
t_{2}+\theta]$ that $$ M(t_{2}+\theta)> M(t_{2})\geq 0,$$ which on
its turn would yield$$M^{2}(t_{2}+\theta)>
M^{2}(t_{2})>\frac{K_{2}}{\varepsilon}.$$ The obtained contradiction
completes the proof of the relation (5.3).

A combination of (5.2) and (5.3) enables us to infer
\begin{equation}
\frac{\sigma}{2}-\varepsilon\leq-\frac{\frac{dM}{dt}}{M^{2}(t)}\leq\frac{\sigma}{2}+\varepsilon
\ \ \ a.e.\ on \ (0,T).
\end{equation}
i.e.
\begin{equation}
\frac{\sigma}{2}-\varepsilon\leq\frac{d}{dt}(\frac{1}{M(t)})\leq\frac{\sigma}{2}+\varepsilon
\ \ \ a.e.\ on \ (0,T).
\end{equation}
For $t\in(t_{2}, T)$ integrating (5.5) on $(t,T)$ to
obtain $$
(\frac{\sigma}{2}-\varepsilon)(T-t)\leq-\frac{1}{M(t)}\leq(\frac{\sigma}{2}+\varepsilon)(T-t),\
t\in (t_{0}, T),$$ that is,$$
\frac{1}{\frac{\sigma}{2}+\varepsilon}\leq-M(t)(T-t)\leq\frac{1}{\frac{\sigma}{2}-\varepsilon},\
t\in(t_{0},T).$$ By the arbitrariness of
$\varepsilon\in(0,-\frac{\sigma}{2})$ the statement of Theorem 5.1
follows.

\bigskip
\noindent\textbf{Acknowledgments} This work was partially supported
by NNSFC (No.11326161, No.11701525).

\end{document}